\documentclass[11pt]{amsart}

\def\comment{}
\def\endcomment{}
\long\def\kill#1\endkill{}

\usepackage{amsmath,amsfonts,amscd,amssymb}
\theoremstyle{plain}
\csname@addtoreset\endcsname{equation}{section}

\newtheorem{theorem}[equation]{Theorem}
\newtheorem{conjecture}[equation]{Conjecture}

\newtheorem{lemma}[equation]{Lemma}
\newtheorem{corollary}[equation]{Corollary}
\theoremstyle{definition}
\newtheorem{remark}[equation]{Remark}

\newtheorem*{notation}{Notation}
\newtheorem{example}[equation]{Example}

\def\X{{\mathcal X}}
\def\Xp{{\mathcal X}_p}
\def\F{{\mathbb F}}

\def\Q{{\mathbb Q}}
\def\Qp{{\mathbb Q}_p}
\def\Zp{{\mathbb Z}_p}

\def\Z{{\mathbb Z}}

\def\triv{{\mathbf 1}}

\def\RC{{\mathcal C}}

\def\smallmatrix#1#2#3#4{
  \genfrac{(}{.}{0pt}{1}{#1}{#3}
  \genfrac{.}{)}{0pt}{1}{#2}{#4}
}

\def\cS{{\mathcal S}}

\def\rksel#1#2#3{\rk_{#3}(#1/#2)}
\def\rkan#1#2{\ord_{s=1}L(#1/#2,s)}
\def\lara{\langle,\rangle}
\def\dlangle{\langle\!\langle}
\def\drangle{\rangle\!\rangle}
\def\llara{\dlangle,\drangle}

\def\neron#1{\omega_{#1}^o}

\newcommand{\vabove}[2]{\genfrac{}{}{0pt}{}{#1}{#2}}

\def\newmathop#1{\expandafter\gdef\csname #1\endcsname{\mathop{\rm #1}\nolimits}}
\newmathop{Tr}
\let\tr\Tr
\newmathop{rk}
\newmathop{ord}
\newmathop{Gal}
\newmathop{Hom}
\newmathop{tors}
\newmathop{coker}
\newmathop{Ind}
\newmathop{Res}
\newmathop{Frob}
\newmathop{GL}
\newmathop{div}
\newmathop{SL}
\newmathop{Aut}
\newmathop{Sel}
\newmathop{End}
\newmathop{Maps}
\newmathop{Re}
\newmathop{Reg}
\def\p{\wp}

\let\lar\longrightarrow
\let\iso\cong
\let\tensor\otimes

\def\leftsemidirect{\mathinner
  {\mathrel\triangleright\joinrel\mathrel{\raise 0.8pt\hbox{$\scriptstyle<$}}}}
\let\leftsemidirect\ltimes
\def\rightsemidirect{\mathinner
  {\mathrel{\raise 0.8pt\hbox{$\scriptstyle>$}}\joinrel\mathrel\triangleleft}}
\let\rightsemidirect\rtimes

\def\dlangle{\langle\!\langle}
\def\drangle{\rangle\!\rangle}

\def\beq{$$\begin{array}{llllllllllllllll}}
\def\eeq{\end{array}$$}

\input cyracc.def       % sha
\font\tencyr=wncyr10
\def\sha{\text{\tencyr\cyracc{Sh}}}

\begin{document}

\let\introdagger\dagger
\title{Self-duality of Selmer groups}
\author{Tim$^\introdagger$ and Vladimir Dokchitser}
\date{June 12, 2008}
\address{Robinson College, Cambridge CB3 9AN, United Kingdom}
\thanks{$^\dagger$Supported by a Royal Society University Research Fellowship}
\email{t.dokchitser@dpmms.cam.ac.uk}
\address{Gonville \& Caius College, Cambridge CB2 1TA, United Kingdom}
\email{v.dokchitser@dpmms.cam.ac.uk}

\comment

\begin{abstract}
%The main result of this paper is the self-duality of Selmer groups:
%if $A$ is an abelian variety over a number field $K$, and
%$F/K$ is a Galois extension with Galois group $G$,
%then the $\Q_pG$-representation naturally associated
%to the $p^\infty$-Selmer group of $A/F$ is self-dual.
%This has several applications to determining parities of Selmer ranks.
%In particular, we prove
%the parity conjecture for Selmer groups of elliptic
%curves $E/\Q$ over abelian extensions of the rationals.

The first part of the paper gives a new proof of self-duality
for Selmer groups:
if $A$ is an abelian variety over a number field $K$, and
$F/K$ is a Galois extension with Galois group $G$,
then the \hbox{$\Q_pG$-}represen\-tation naturally associated
to the $p^\infty$-Selmer group of $A/F$ is self-dual.
The second part describes a method for obtaining information
about parities of Selmer ranks from the local Tamagawa numbers of $A$
in \hbox{intermediate} extensions of $F/K$.
\end{abstract}

%\llap{.\hskip 10cm} \vskip -0.4cm
\maketitle
\def\introdagger{{}}

%\tableofcontents

\section{Introduction}

Let $F/K$ be a Galois extension of number fields with Galois group $G$,
and $A$ an abelian variety defined over $K$. The
action of $G$ on the $F$-rational points of $A$
defines a $\Q G$-representation $A(F)\tensor\Q$, which
in particular recovers the Mordell-Weil ranks of $A$ over all intermediate
extensions.

Now if $p$ is a prime number, there is an analogous picture for Selmer groups.
Let
$$
  \small
  \X=\Xp(A/F)=(\text{Pontryagin dual of the $p^\infty$-Selmer group of $A/F$})\tensor\Q_p.
$$
It is a $\Q_p$-vector space whose dimension is the $p^\infty$-Selmer rank
$\rksel AFp$,
the Mordell-Weil rank plus the number of copies of $\Q_p/\Z_p$ in
the Tate-Shafarevich group $\sha(A/F)$.
Moreover, it is a $G$-representation, and it also recovers
$\rksel ALp$ for intermediate extensions as
the dimension of $\X^{\Gal(F/L)}$.
In view of the conjectural finiteness of the
Tate-Shafarevich group,
$\X$ and $A(F)\tensor\Q_p$ should be isomorphic
as $G$-representations.

Decompose
$
  \X \iso \bigoplus_{\rho} \rho^{\oplus m_\rho}
$
into $\Q_p G$-irreducible constituents.
The structure of $\X$ is encoded in
the multiplicities $m_\rho$.
Our main result (Theorem \ref{imrho} below) may be stated as
$$
  \sum_{\rho\in S_\Theta} m_\rho \mod 2 \quad=\quad \text{(explicit local data)}
$$
for certain special sets $S_\Theta$ of representations.

A result of this type is proved in \cite{TV-S} for Mordell-Weil groups
under the assumption that $\sha(A/F)$ is finite, and it relies crucially
on the existence of the height pairing.
The Selmer group also admits a non-degenerate \hbox{$G$-invariant} pairing
$$
  \langle,\rangle:\>\> \X \times \X \lar \Qp.
$$
Equivalently,
\begin{theorem}
\label{iselfdual}
$\X$ is self-dual as a $\Qp G$-representation.
\end{theorem}

%This is essentially a consequence of Poitou-Tate duality,
%but as it is difficult to find in the literature in this form
%we present a new proof of it.

\noindent
This is essentially a consequence of Poitou-Tate duality,
and may be deduced
using the methods of Greenberg, see \cite{Gre2} Prop. 2.
A proof in the general context of Bloch-Kato Selmer groups is given by
Nekov\'a\v r in \cite{Nek} 12.5.9.5(iv).
We present an alternative proof in \S2.

Let us briefly record two straightforward consequences
of self-duality, whose proofs we postpone to \S 2.
The first of these relies on the parity conjecture for elliptic curves
over $\Q$, which is known thanks to the work of
Birch--Stephens~\cite{BS}, Greenberg~\cite{Gre} and Guo~\cite{Guo}~($E$ CM),
Monsky \cite{Mon} ($p=2$), Nekov\'a\v r~\cite{Nek} ($p$ potentially ordinary),
Kim~\cite{Kim} ($p$ supersingular) and \cite{TV-S} ($p$~odd).
In particular, Theorem \ref{iabext} is proved in \cite{Nek}
when $E$ has potentially ordinary reduction at $p$.

\begin{theorem}[$=\,$Theorem \ref{abext}]
\label{iabext}
Let $E/\Q$ be an elliptic curve. For every abelian extension $F/\Q$
and every prime $p$,
$$
  \rksel EFp \equiv \rkan EF \mod 2.
$$
\end{theorem}

\begin{theorem}[$=\,$Corollary \ref{corodd}]
Let $A/K$ be an abelian variety, and suppose
$F/K$ is Galois of odd degree. Then $\rksel AFp\equiv\rksel AKp\mod 2$.
\end{theorem}

\medskip

Returning % down to Earth
to the representation-theoretic structure of $\X$,
suppose we are given a relation between induced representations (fixed for
the rest of the introduction),
$$
  \Theta:\>\>  \bigoplus_i\Ind_{H_i}^G\triv\iso\bigoplus_j\Ind_{H'_j}^G\triv\qquad (H_i, H'_j<G).
$$
Observe that Artin formalism forces an equality of $L$-functions
$$
  \prod L(A/F^{H_i},s) = \prod L(A/F^{H'_j},s),
$$
and that the conjectural Birch--Swinnerton-Dyer formula at $s=1$ implies a relation
between the arithmetic invariants of $A$ over these fields.
As explained in \cite{TV-S} \S2.2,
most of these cancel modulo rational squares, leading~to

\begin{conjecture}[$\square$-Conjecture]
Suppose $A$ has a principal polarisation induced by a $K$-rational divisor
(e.g. $A$ is an elliptic curve.) Then
$$
  \frac{\prod_i \Reg(A/F^{H_i})}{\prod_j \Reg(A/F^{H'_j})}
  \equiv
  \frac{\prod_i C(A/F^{H_i})}{\prod_j C(A/F^{H'_j})}
  \mod
  \Q^{*2}.
$$
\end{conjecture}

\noindent
Here $\Reg$ is the regulator, and $C$ is a product of local terms
(essentially Tamagawa numbers, see the list of notation below).
The assumption on $A$ guarantees that
$\sha$ modulo its divisible part has square order.

We showed (\cite{TV-S} Cor. 2.5) that this conjecture is implied
by the Shafarevich-Tate conjecture on the finiteness of $\sha$.
Now we use the pairing $\lara$ on $\X$ to prove the
corresponding unconjectural statement for Selmer groups.
For every subgroup $H$ of $G$ define
\def\Regp{\Reg_p^{\scriptscriptstyle\lara}}
$$
  \Regp(A/F^H) =
    \det\bigl(\tfrac{1}{|H|}\langle,\rangle\bigm|\X^H\bigr)\>\> \in \Q_p^*/\Q_p^{*2},
$$
computing the determinant on any basis of $\X^H$.
This, including the scaling, is analogous to the definition
of $\Reg(A/F^H)$ as the determinant of the
height pairing
on $A(F^H)$ modulo torsion.

\begin{theorem}[$=$Theorem \ref{selcompat}]
\label{iselcompat}
Suppose $A/K$ is principally polarised.
If $p\!=\!2$, assume furthermore that the principal polarisation on $A$
is induced by a $K$-rational divisor. Then
$$
  \ord_p
  \frac{\prod_i \Regp(A/F^{H_i})}{\prod_j \Regp(A/F^{H'_j})}
  \equiv
  \ord_p
  \frac{\prod_i C(A/F^{H_i})}{\prod_j C(A/F^{H'_j})}
  \mod 2.
$$
\end{theorem}

The theorem may be used to express parities of Selmer ranks in terms of local
invariants. The crucial point is that the left-hand side
depends only on $\X$
as a $\Q_pG$-representation, and not on the pairing $\lara$.
Let $\cS$ be the set of self-dual $\Qp G$-representations, which
are either irreducible or of the form $T\oplus T^*$ for some irreducible
$T\not\iso T^*$ ($T^*$ is the contragredient of $T$).
Any self-dual $\Qp G$-representation can be uniquely
decomposed into such constituents.
For every $\rho\in\cS$, pick a non-degenerate $G$-invariant pairing
$\llara$ on $\rho$, and define the {\em regulator constant\/}
$$
  \RC(\Theta,\rho)=\frac
     {\prod_i\det(\tfrac{1}{|H_i|}\llara|\rho^{H_i})}
     {\prod_j\det(\tfrac{1}{|H'_j|}\llara|\rho^{H'_j})}
  \in \Q_p^*/\Q_p^{*2}.
$$
%(This does not depend on the pairing by Theorem \ref{regconst}.)
Consider the set
$$
  \cS_\Theta = \{\rho\in\cS \,|\, \ord_p\RC(\Theta,\rho)\equiv 1\mod 2\}.
$$
This is a ``computable combination of representations''
in the following sense:

\begin{theorem}
\label{imrho}
Suppose $A/K$ is principally polarised, and if $p\!=\!2$ assume furthermore
that the polarisation is induced by a $K$-rational divisor.
There is a decomposition $\X \iso \bigoplus_{\rho\in\cS} \rho^{\oplus m_\rho}$,
and
$$
  \sum_{\rho\in\cS_\Theta} m_\rho \equiv
  \ord_p
  \frac{\prod_i C(A/F^{H_i})}{\prod_j C(A/F^{H'_j})}
  \mod 2.
$$
\end{theorem}

\begin{proof}
By self-duality, such a decomposition exists.
Take the obvious pairing $\lara$ on $\X$ coming from $\llara$
on its constituents and apply Theorem \ref{iselcompat}.
\end{proof}

Note that in practice the right-hand side is very explicit:
it can be computed for elliptic curves by Tate's algorithm (\cite{Sil2} IV.9),
for semistable abelian varieties from the monodromy pairing (\cite{Gro} \S10),
and for Jacobians of curves using the intersection pairing (\cite{BL} \S1).

In Examples \ref{ex1}--\ref{ex3}
we will illustrate this theorem for specific relations when
$G=D_{2p}, \hbox{$C_p\!\rightsemidirect\!C_{p-1}$}$ and $\GL_2(\F_p)$.
The first two were already considered in \cite{TV-S}, but
required ad hoc constructions of isogenies in absence of Theorem \ref{imrho},
see \cite{TV-S} Thm. 4.11, Prop. 4.17.
The $D_{2p}$-extensions were also studied by
Mazur and Rubin, who give another local expression for the same
parity of Selmer ranks, see \cite{MR} Thm.\ A.
We do not have an intrinsic description of the sets $\cS_\Theta$ for a
general group $G$.

In the forthcoming paper \cite{TV-T} we will expand on the properties of
permutation relations and regulator constants, and address
the question of the compatibility of Theorem \ref{imrho} with root numbers.
Under reasonably mild hypotheses on $A$
(e.g. when $A$ is an elliptic curve whose additive primes above 2 and 3 are
unramified in $F/K$),
we will
show that the parity of $\sum_{\rho\in\cS_\Theta} m_\rho$ is
indeed determined by the (conjectural) sign in the functional equation of a
corresponding $L$-function, as predicted by the parity conjecture.

\medskip

\noindent {\bf Acknowledgements.\/}
We would like to thank Christian Wuthrich for his comments.

\begin{notation}
Throughout the paper we fix

\noindent
\begin{tabular}{lll}
&$F/K$                 & Galois extension of number fields \cr
&$G$                   & $\Gal(F/K)$                       \cr
&$A/K$                 & abelian variety with a fixed regular non-zero exterior form $\omega$ \cr
&$p$                   & prime number                      \cr
\end{tabular}

%\pagebreak
\noindent
For an intermediate field $K\subset L\subset F$,
we use the following notation:

\noindent
\begin{tabular}{lll}
&$\Sel_{p^\infty}(A/L)$ & the $p^\infty$-Selmer group $\varinjlim\Sel_{p^n}(A/L)$\cr
&$X_p(A/L)$            & Pontryagin dual $\Hom(\Sel_{p^\infty}(A/L),\Q_p/\Z_p)$ modulo torsion\cr
&$\Xp(A/L)$            & $X_p(A/L)\tensor_{\Z_p}\Q_p$\cr%, viewed as a $\Qp G$-representation\cr
&$\rksel ALp$          & $p^{\infty}$-Selmer rank of $A/L$, i.e. $\dim \Xp(A/L)$ \cr
%&                      & $\rk A(L)$ plus the $\Z_p$-corank of $\sha(A/L)$ \cr
&$C(A/L)$              & $\prod c_v |{\omega}/{\neron{v}}|_v$, where
   the product is taken over all primes \cr&&  of $L$, $c_v$ is the local
   Tamagawa number, $\neron{v}$ the N\'eron \cr&& differential and
   $|\cdot|_v$ the normalised absolute value.
\end{tabular}

\noindent
By convention, permutation modules $\Q[G/H]\iso\Ind_H^G\triv$ come
with a standard basis of elements of $G/H$.
With respect to this basis, the identity matrix defines
a $G$-invariant pairing.
\end{notation}

%\pagebreak[3]

\section{Self-duality}

The purpose of this section is to prove

\begin{theorem}[$=$Theorem \ref{iselfdual}]
\label{selfdual}
Suppose $F/K$ be a Galois extension of number fields with Galois group $G$,
$A/K$ an abelian variety, and $p$ a prime number.
Then $\Xp(A/F)$ is self-dual as a $\Qp G$-representation.
\end{theorem}

To begin with,
if $M$ is a finitely generated $\Z G$-module, there is a naturally associated
abelian variety $A\tensor M$ over $K$ (\cite{MilO} \S2).
When $L/K$ is an intermediate extension, $A\tensor\Z[G/\Gal(F/L)]$ is
the Weil restriction of scalars $W_{L/K}A$.
An injection of $G$-modules $\phi: M_1\to M_2$ with finite cokernel
induces a $K$-isogeny $f_\phi: A\tensor M_1\to A\tensor M_2$.
If $A$ is principally polarised and $M_1$ and $M_2$ are permutation modules,
then $A\tensor M_i$ carry induced polarisations, with respect to which
the dual isogeny $f_\phi^t: (A\tensor M_2)^t\to (A\tensor M_1)^t$
comes from the transposed matrix $\phi^t$ in the standard bases (\cite{TV-S} \S4.2).
In other words, $f_{\phi^t}=(f_\phi)^t$.

An isogeny $f: A\to B$ induces a $G$-invariant map $X_p(B/K)\to X_p(A/K)$,
which is an isomorphism when tensored with $\Qp$.
Following \cite{TV-S} \S4, write
\beq
  Q(f) = & |\coker(f: A(K)/A(K)_{\tors} \to B(K)/B(K)_{\tors})|\>\>\times\cr
            &  \>\>\times\>\>
           |\ker(f: \sha(A)_{\div}\to \sha(B)_{\div})|\>,\cr
\eeq
where $\sha_{\div}$ denotes the divisible part of $\sha$.
Then $Q(f)$ is multiplicative under composition of isogenies, and
its $p$-part is the size of the cokernel of $f: X_p(B/K)\to X_p(A/K)$.

Recall the Selmer group analogue of the invariance of the
Birch--Swinnerton-Dyer quotient under isogenies:

\begin{theorem}[\cite{TV-S} Thm. 4.3]
\label{thmsel}
Let $A, B/K$ be abelian varieties given with regular non-zero exterior
forms $\omega_A, \omega_B$, and suppose $\phi: A\to B$ is a \hbox{$K$-isogeny}.
Writing $\sha_0(A/K)$ for $\sha(A/K)$
modulo its divisible part and
$$
  \Omega_A=
  \prod\limits_{\vabove{v|\infty}{\text{real}}}
  {\int_{A(K_v)} |\omega_A|}\cdot
  \prod\limits_{\vabove{v|\infty}{\text{complex}}}
  {\int_{A(K_v)}\! \omega_A\!\wedge\! \bar\omega_A}
$$
and similarly for $B$, we have
\begingroup
$$
  \frac{|B(K)_{\tors}|}{|A(K)_{\tors}|}
  \frac{|B^t(K)_{\tors}|}{|A^t(K)_{\tors}|}
  \frac{C(A/K)}{C(B/K)}
  \frac{\Omega_A}{\Omega_B}
  \prod_{l|{\deg}\phi}\!
    \frac{|\sha_0(A)[l^\infty]|}{|\sha_0(B)[l^\infty]|}
     \!=\!
  \frac{Q(\phi^t)}{Q(\phi)}.
$$
\endgroup
\end{theorem}

To prove Theorem \ref{selfdual}, we establish the analogous
statement for $\Xp(A\tensor\Z[G]/K)$ in place of $\Xp(A/F)$, and then
show that the two are isomorphic.

\begin{theorem}
\label{seldual}
Let $F/K$ be a Galois extension of number fields,
$A/K$ an abelian variety and $p$ a prime number.
Then $\Xp(A\tensor\Z[G])$ is a self-dual $\Gal(F/K)$-representation.
%$S=\Hom_{\Qp}(X_p(A/F)\tensor_{\Zp}\Qp,\Qp)$ the double dual Selmer
\end{theorem}

\begin{proof}
Let $M=\Z[G]$ and $\X=\Xp(A\tensor M/K)$.
The idea is that for a self-isogeny
$$
  f:\> A\tensor M \lar A\tensor M
$$
we have $Q(f)=Q(f^t)$ by Theorem \ref{thmsel}. We will construct an $f$ whose
$Q$ recovers the multiplicity of a given representation in
$\X$, and $Q(f^t)$ recovers the multiplicity of its dual.

To begin with, after passing to an isogenous abelian variety if necessary,
we may assume $A$ that is principally polarised.
We have to show that for every $g\in G$, its eigenvalues
on $\X$ come in inverse pairs. Restricting to the
subgroup generated by $g$, we will also assume that $G$ is cyclic.

Let $\tau_i$ be the distinct $\Qp$-irreducible representations of $G$
and write \linebreak $m_{\tau_i}(\X)$ for their multiplicity in $\X$.
Choose $G$-invariant $\Zp$-sublattices $\Lambda_{\tau_i}$ of
\hbox{$X_p(A\tensor M/K)$}~with
$\Lambda_{\tau_i}\tensor_{\Zp}\Qp\iso\tau_i{}^{\oplus m_{\tau_i}(\X)}$,
so $\bigoplus \Lambda_{\tau_i}$ is of finite index.

For $\tau\in\{\tau_i\}$,
let $P_\tau = \sum_{g\in G} \tr (\tau (g)) g^{-1} \in \Zp[G]$
be ``$|G|$ times the projector onto the $\tau$ component'' operator.
%? It is invertible in the group algebra, e.g. use Cayley-Hamilton
Elements in an open neighbourhood $U$ of
$|G| + (p-1)P_\tau$ act as $|G|$ times an isomorphism on all $\Lambda_{\tau_i}$
save $\Lambda_\tau$,
where they act by $p|G|$ times an isomorphism. Similarly there is such
a neighbourhood $U^*$ for $\tau^*$.
%``Unit'' means $\Zp$-integral matrix in our basis,
%sufficiently close to identity (and, in particular, prime-to-$p$ determinant).
Since the $\Qp$-linear map on $\Qp[G]$ determined by $g\to g^{-1}$ for
$g\in G$ is continuous and sends $P_\tau$ to $P_{\tau^*}$, we
can choose $\Phi=\sum_g x_g g \in\Z[G]\cap U$
with $\Phi^*=\sum_g x_g g^{-1}\in U^*$.

Since $G$ is commutative,
$\Phi$ defines a $G$-endomorphism $\phi$ of $M$. Considering
its action on the $\Lambda_{\tau_i}$,
$$
  Q(f_\phi) = Q(f_{|G|}) \> p^{m_\tau(\X)\dim\tau}.
$$
%(from the definition of $Q$ --- integral units have trivial cokernels).
%using that $\Z$-integral matrices with $p$-unit determinant have $Q=1$.
Also, $\Phi^*$ corresponds to $\phi^t$ in $\End(M)$ with respect to
the standard basis of $M=\Z[G]$, so
$$
  Q(f_{\phi^t}) = Q(f_{|G|}) \> p^{m_{\tau^*}(\X)\dim\tau^*}.
$$
%We have to show that its multiplicity $m_\tau(\X)$ in $\X$
%coincides with that of $\tau^*$.
Since $Q(f_\phi)/Q(f_{\phi^t})=1$ by Theorem \ref{thmsel},
$m_\tau(\X)=m_{\tau^*}(\X)$. As this holds for all $\tau$,
the asserted self-duality follows.
\end{proof}

\begin{lemma}%[Equivariant Shapiro]
\label{equishapiro}
Let $A/K$ be an abelian variety, and
$W$ the Weil restriction $W_{F/K}A$.
There is a canonical isomorphism of $G$-modules
$$
  \Sel_{p^n}(W/K) \lar \Sel_{p^n}(A/F),
$$
where $G$ acts on $\Sel_{p^n}(W/K)$ via automorphisms of $W/K$
and on $\Sel_{p^n}(A/F)$ by its usual action on $H^1(F,A[p^n])$.
In particular, $\Xp(W/K)\iso\Xp(A/F)$ as $G$-representations.
\end{lemma}

\begin{proof}
By Milne \cite{MilA} p.178 (a), we have $W[p^n]=\Ind_K^F A[p^n]$, and
the isomorphism given by Shapiro's lemma,
$$
  %\alpha:
  H^1(K,\Ind_K^F A[p^n]) \lar H^1(F,A[p^n]),
$$
descends to an isomorphism of Selmer groups. It is easy to check that it
is compatible with the $G$-action.
\kill
We consider $\Ind_K^F A[p^n]$ as functions $f: G \to A[p^n]$,
with $\Gal(\bar K/K)$-action $(\gamma\cdot f)(x)=\gamma\cdot f(x \bar\gamma)$
and a commuting $G$-action $(g\cdot f)(x)=f(g^{-1}\cdot x)$.

Take an element of $H^1(K,\Ind_K^F A[p^n])$, represented by
a cocycle $c: \Gal(\bar K/K)\to \Ind_K^F A[p^n]$
mapping $\sigma$ to $c_\sigma \in \Maps(G,A[p^n])$.
For $g\in G$,
$$
  (g\cdot c)_\sigma = (x \mapsto c_\sigma(g^{-1}x)).
$$
For the action of $g\in G$ on $H^1(F,A[p^n])$,
with $\gamma\in \Gal(\bar K/K)$
a lift of $g$,
$$
  (g\cdot d)(\tau) = \gamma \cdot d(g^{-1} \tau g).
$$
Then $\alpha(c)$ is given by $\sigma \mapsto c_\sigma(1)$.

We have
$$
  g\cdot\alpha(c): \>\>\tau \longmapsto \gamma \cdot (\sigma\mapsto c_\sigma(1))(g^{-1} \tau g)
    = \gamma\cdot c_{g^{-1} \tau g}(1)
$$
and
$$
  \alpha(g\cdot c):\>\>\tau \longmapsto
    (g\cdot c)_\tau(1) = (x \mapsto c_\tau(g^{-1}x))(1) = c_\tau(g^{-1}).
$$
But
$$
\begin{array}{lll}
  \gamma\cdot c_{g^{-1} \tau g}(1) &=&
    gg^{-1}\cdot c_{\tau g}(1\cdot g^{-1}) + g\cdot c_{g^{-1}}(1)\cr
    &=&\tau\cdot c_g(g\cdot\bar\tau) + c_\tau(g^{-1}) + g\cdot c_{g^{-1}}(1)\cr
    &=&\tau\cdot c_g(g) + c_\tau(g^{-1}) + g\cdot c_{g^{-1}}(1)\cr
    &\equiv&c_g(g) + c_\tau(g^{-1}) + g\cdot c_{g^{-1}}(1)\cr
    &\equiv&c_\tau(g^{-1}),\cr
\end{array}
$$
using the cocycle condition in the first two equalities,
$\bar\tau=1$, and working modulo coboundaries in the last two.
\endkill
\end{proof}

This completes the proof of Theorem \ref{selfdual}.

%\newpage

\begin{corollary}
\label{corodd}
Let $A/K$ be an abelian variety, and suppose
$F/K$ is Galois of odd degree. Then $\rksel AFp\equiv\rksel AKp\mod 2$.
\end{corollary}

\begin{proof}
Since $G=\Gal(F/K)$ has odd degree, its only self-dual irreducible
representation is the trivial one. (Their number coincides with the number
of self-inverse conjugacy classes of $G$, but these have odd order and,
except for the trivial class, have no self-inverse elements.)
So $\Xp(A/F)\iso \Xp(A/K)\oplus$(even-dimensional representation).
\end{proof}

\begin{remark}
The $p$-parity conjecture for abelian varieties asserts that
$(-1)^{\dim\Xp(A/K)}$ coincides with the root number $w(A/K)$, the
(conjectural) sign in the functional equation for $L(A/K,s)$.
If $F/K$ is Galois of odd degree, then $w(A/F)=w(A/K)$
(see \cite{TatN} 3.4.7, 4.2.4), so the $p$-parity conjecture holds for $A/F$
if and only if it holds for $A/K$, by the corollary.
\end{remark}

\begin{corollary}
\label{corabelian}
Let $A/K$ be an abelian variety, and suppose $F/K$ is abelian.
Let $K(\sqrt{d_i})$ be the quadratic extensions of $K$ in $F$ and write
$A_i$ for the corresponding quadratic twists of $A$. Then
$$
  \rksel AFp\equiv\rksel AKp+{\textstyle\sum_i\,} \rksel{A_i}Kp \mod 2.
$$
\end{corollary}

\begin{proof}
The self-dual irreducible complex representations
of $\Gal(F/K)$ are characters of order 1 or 2.
\end{proof}

\begin{theorem}[Parity Conjecture in abelian extensions]
\label{abext}
Let $E/\Q$ be an elliptic curve. For every abelian extension $F/\Q$
and every prime $p$,
$$
  \rksel EFp \equiv \rkan EF \mod 2.
$$
\end{theorem}

\begin{proof}
With the notation from Corollary \ref{corabelian},
%The self-dual irreducible complex representations
%of $\Gal(F/K)$ have order 1 or 2. Write $E_i$ for the twists of $E$ by
%these characters.
%nd they correspond to quadratic fields
%Let $\chi_i: \Gal(K_i/\Q)\to\{\pm 1\}$ be the non-trivial characters
%of the all quadratic fields inside $F$, and $E_i$ the twists of $E$ by
%$\chi_i$.
%By the self-duality of $\Xp(E/F)$,
$$
  \rksel EFp\equiv\rksel EKp+{\textstyle\sum_i\,} \rksel{E_i}Kp \mod 2.
$$
In view of the Parity Conjecture for elliptic curves over $\Q$
(\cite{TV-S} Theorem 1.4), the right-hand side agrees with
the corresponding analytic ranks. By the functional equation,
$\ord_{s=1}L(E,\tau,s)=\ord_{s=1}L(E,\tau^*,s)$ for every
character $\tau$ of $\Gal(F/K)$, so
$$
  \rkan EF \equiv \rkan{E}\Q + {\textstyle\sum_i} \rkan{E_i}\Q \mod 2.
$$
\end{proof}

\section{Regulator constants for Selmer groups}

\def\Regpf#1{\det\bigl(\tfrac{1}{|#1|}\langle,\rangle\bigm|{\X}^{#1}\bigr)}

The central result of this section is

\begin{theorem}[$=$Theorem \ref{iselcompat}]
\label{selcompat}
Let $G\!=\!\Gal(F/K)$ and let $p$ be a prime number. Suppose
$A/K$ is principally polarised, and if $p\!=\!2$,
assume furthermore that the principal polarisation on $A$
is induced by a $K$-rational divisor.
If $H_i,H'_j<G$ satisfy
$\oplus\Ind_{H_i}^G\triv\iso\oplus\Ind_{H'_j}^G\triv$, then
for every
non-degenerate $G$-invariant $\Q_p$-bilinear pairing $\lara$ on $\X=\Xp(A/F)$,
$$
  \ord_p
  \frac{\prod_i \Regpf{H_i}}{\prod_j \Regpf{H'_j}}
  \equiv
  \ord_p
  \frac{\prod_i C(A/F^{H_i})}{\prod_j C(A/F^{H'_j})}
  \mod 2.
$$
\end{theorem}
\begin{proof}
Write $S_1=\coprod_i G/H_i$ and $S_2=\coprod_j G/H'_j$.
Since $\Q[S_1]\iso\Q[S_2]$, there is a $G$-injection
$\phi: \Z[S_1]\to\Z[S_2]$ with
finite cokernel, and it induces maps
%\newpage\noindent
on abelian varieties
$$
  f_\phi: A\tensor\Z[S_1] \lar A\tensor\Z[S_2], \qquad
  f_{\phi^t}: A\tensor\Z[S_2] \lar A\tensor\Z[S_1].
$$
Applying Theorem \ref{thmsel} modulo rational squares
(see also \cite{TV-S} Cor. 4.5),
$$
\ord_p \frac{\prod_i C(A/F^{H_i})}{\prod_j C(A/F^{H'_j})}
  \equiv \ord_p Q(f_{\phi\phi^t})\mod 2.
$$
It remains to justify the last two steps in the following chain of equalities:
$$
\begin{array}{rcl}
  \ord_p Q(f_{\phi\phi^t})
  &=&
  \ord_p\coker(f_{\phi\phi^t}\,|\,X_p(A\tensor\Z[S_2])) \cr
  &=&
  \ord_p\det(f_{\phi\phi^t}\,|\,\Xp(A\tensor\Z[S_2])) \cr
  &{\buildrel Cor. \ref{corselcd}\over=}&
  \ord_p\det((\phi\phi^t)^*\,|\, \Hom_G(\Z[S_2],\X)  ) \cr
  &{\buildrel Lem. \ref{qvsreg}\over\equiv}&
 % \displaystyle
  \ord_p\frac
  {\prod_i \det(\frac{1}{|H_i|}\langle,\rangle|\X^{H_i})}
  {\prod_j \det(\frac{1}{|H'_j|}\langle,\rangle|\X^{H'_j})}
  \mod 2.
\end{array}
$$
\vskip -8pt
\end{proof}

\begin{lemma}
\label{qvsreg}
Suppose $S_1, S_2$ are finite $G$-sets, and $\phi: \Z[S_1]\to\Z[S_2]$
is a $G$-injection with finite cokernel. Write $\phi^t: \Z[S_2]\to\Z[S_1]$
for its transpose in the standard basis.
For a $\Qp G$-representation $V$, write
$$
  \phi^*: \Maps_G(S_2,V) \lar \Maps_G(S_1,V)
$$
for the pullback of maps, and similarly for $\phi^t$ and $\phi\phi^t$.

If $V$ has a $G$-invariant $\Qp$-bilinear non-degenerate pairing $\lara$,
and we write $S_1=\coprod G/H_i, S_2=\coprod G/H'_j$ with $H_i,H'_j<G$, then
$$
  \ord_p\det ((\phi\phi^t)^*) \equiv \ord_p
  \frac{\prod_i \det\bigl(\tfrac{1}{|H_i|}\langle,\rangle\bigm|V^{H_i}\bigr)}%
       {\prod_j \det\bigl(\tfrac{1}{|H'_j|}\langle,\rangle\bigm|V^{H'_j}\bigr)}
  \mod 2.
$$
\end{lemma}

\begin{proof}
We may identify
$$
   V^H =\Maps_G(G/H,V) =
     \{f: G/H\to V \>|\> f(g\cdot s)=g\cdot f(s)\}
$$
with the map from right to left given by $f\mapsto f(1)$.

If $S$ is a finite $G$-set, define an inner product on functions $S\to V$ by
$$
  (f_1,f_2) = \tfrac{1}{|G|} \sum_{s\in S} \langle f_1(s), f_2(s) \rangle.
$$
For $S=G/H$
it agrees with the inner product $\tfrac{1}{|H|}\langle,\rangle$ on $V^H$ via
the identification above. Indeed, if $f_1(1)=v_1$ and $f_2(1)=v_2$ then
$$
\begin{array}{rcl}
  (f_1,f_2) &=&
            \displaystyle
            \tfrac{1}{|G|} \sum_{s\in G/H} \langle f_1(s), f_2(s) \rangle
            = \tfrac{1}{|G|} \sum_{s\in G/H} \langle s\cdot f_1(1), s\cdot f_2(1) \rangle\cr
            &=&
            \displaystyle
            \tfrac{1}{|G|} \sum_{s\in G/H} \langle v_1,v_2 \rangle
            = \tfrac{1}{|H|} \langle v_1,v_2 \rangle.\cr
\end{array}
$$
Therefore for $S_1=\coprod G/H_i$ (and similarly for $S_2$),
$$
  \ord_p
  \prod_i \det\bigl(\tfrac{1}{|H_i|}\langle,\rangle\bigm|V^{H_i}\bigr)
  \equiv
  \ord_p
  \det\bigl(\, (,) \bigm| \Maps_G(S,V) \,\bigr)
  \mod 2.
$$
An elementary computation shows that $\phi^*$ and $(\phi^t)^*$ are
adjoint with respect to $(,)$, so picking a basis $\{e_k\}$ of $\Maps_G(S_2,V)$,
\beq
  \ord_p\det\bigl(\, (,) \bigm| \Maps_G(S_1,V) \,\bigr)
  \equiv
  \ord_p
  \det\bigl((\phi^*e_k,\phi^*e_l)_{k,l}\bigr) \cr
  \equiv
  \ord_p\>
  \det\bigl((e_k,(\phi^t)^*\phi^*e_l)_{k,l}\bigr) \cr
  \equiv
  \ord_p\>
  \det((\phi\phi^t)^*)\det\bigl(\, (,) \bigm| \Maps_G(S_2,V) \,\bigr)
  \mod 2.
\eeq
\end{proof}

\begin{lemma}
\label{lemselcd}
Let $A/K$, $F/K$ and $G$ be as in Theorem \ref{selcompat}, and
suppose $\phi: \Z[S_1]\to\Z[S_2]$ is an injection of $G$-permutation modules
with finite cokernel. There are natural vertical maps with finite kernels and
cokernels that make the diagram
$$
  \begin{CD}
  \Hom_G(\Z[S_1],\Sel_{p^\infty}(A/F)) & @>{\phi^*}>> & \Hom_G(\Z[S_2],\Sel_{p^\infty}(A/F)) \\
  @AAA && @AAA \\
  \Sel_{p^\infty}(A\tensor \Z[S_1]/K) & @>{f_\phi}>> & \Sel_{p^\infty}(A\tensor \Z[S_2]/K) \\
  \end{CD}
$$
commute. Here $\phi^*$ is the pullback of maps induced by $\phi$.
\end{lemma}
\begin{proof}
Writing $S_1\!=\!\coprod G/H_i$, $S_2\!=\!\coprod G/H'_j$, we
have a commutative diagram
$$
  \begin{CD}
  \Hom_G(\prod\Z[ G/H_i],H^1(F,A[p^n])) & @>{\phi^*}>> & \Hom_G(\prod\Z[ G/H'_j],H^1(F,A[p^n])) \\
  @V{{\rm eval}(1,1,...,1)}V\iso V && @V{{\rm eval}(1,1,...,1)}V\iso V \\
  \prod H^1(F,A[p^n])^{H_i} &&&& \prod H^1(F,A[p^n])^{H'_j} \\
  @A{\Res}AA && @A{\Res}AA \\
  \prod H^1(F^{H_i},A[p^n]) &&&& \prod H^1(F^{H'_j},A[p^n]) \\
  @A{{\rm eval}(1,1,...,1)}A\iso A && @A{{\rm eval}(1,1,...,1)}A\iso A \\
  H^1(K,\Hom(\prod\Z[ G/H_i],A[p^n])) & @>{\phi^*}>> & H^1(K,\Hom(\prod\Z[ G/H'_j],A[p^n]))\\
  \end{CD}
$$
The restriction maps $\Res$ have bounded kernels
and cokernels with respect to $n$.
Taking the limit (and using a similar local diagram)
proves the lemma.
\end{proof}

\begin{corollary}
\label{corselcd}
In the situation of Lemma \ref{lemselcd}, there is a commutative \hbox{diagram}
$$
  \begin{CD}
  \Hom_G(\Z[S_1],\Xp(A/F)) & @<{(\phi^t)^*}<< & \Hom_G(\Z[S_2],\Xp(A/F)) \\
  @V{\iso}VV && @V{\iso}VV \\
  \Xp(A\tensor \Z[S_1]/K) & @<{f_\phi}<< & \Xp(A\tensor \Z[S_2]/K) \\
  \end{CD}
$$
\end{corollary}

\begin{proof}
Apply Pontryagin duals $\Hom_{\Zp}(\cdot,\Qp/\Zp)$ to the diagram
of Lemma \ref{lemselcd}, and tensor with $\Qp$.
The fact that the Pontryagin duals in the top row are what the corollary
asserts they are is general nonsense:
for rings $R, S$ and modules
${\mathcal A}_R, {}_R{\mathcal B}_S, {\mathcal C}_S$,
there is a natural isomorphism (see e.g. \cite{Mac} Thm. V.3.1, p.144)
$$
  \eta: \Hom_S({\mathcal A}\tensor_R {\mathcal B},{\mathcal C})\iso\Hom_R({\mathcal A},\Hom_S({\mathcal B},{\mathcal C}))
$$
of abelian groups, defined for $h: {\mathcal A}\tensor_R {\mathcal B}\to {\mathcal C}$ by
$[(\eta h)a](b)=h(a\tensor b)$. Apply this with
$R=\Z[G]$, $S=\Zp$,
${\mathcal A}=\Hom_\Z(\Z[S_i],\Z)$, ${\mathcal B}=\Sel_{p^\infty}({\mathcal A}\tensor \Z[S_i]/K)$ and ${\mathcal C}=\Qp/\Zp$.
Note that $\phi: \Z[S_1]\to\Z[S_2]$ induces the transpose map
$$
  \phi^t: \Z[S_2]=\Hom_\Z(\Z[S_2],\Z)\to \Hom_\Z(\Z[S_1],\Z)=\Z[S_1].
$$
\end{proof}

We illustrate the applications of Theorem \ref{imrho} with a few examples.
As in the theorem, suppose $A/K$ is a principally polarised abelian variety,
$F/K$ a Galois extension, and $p$ a fixed prime.
Write $m_\tau$ for the multiplicity of $\tau$ in $\Xp(A/F)$,
and $\Theta=\sum H_i-\sum H'_j$ for a relation
$\Theta: \bigoplus_i\Ind_{H_i}^G\triv\iso\bigoplus_j\Ind_{H'_j}^G\triv$.

\begin{example}
\label{ex1}
Suppose $p$ is an odd prime and $G=\Gal(F/K)\iso D_{2p}$ is dihedral.
Let $M, L$ be intermediate fields of degree 2 and $p$ over $K$, respectively.
The group $G$ has three $\Q_p$-irreducible representations,
which are all self-dual:
trivial $\triv{}$, sign $\epsilon$ and $(p\!-\!1)$-dimensional $\rho$.
There is a unique (up to multiples) relation between permutation representations,
$$
  \Theta=\{1\} - \>2\Gal(F/L) - \Gal(F/M) + 2\>G.
$$
A simple computation shows that
$\RC(\Theta,\triv)=\RC(\Theta,\epsilon)=\RC(\Theta,\rho)=p$, so
$\cS_\Theta=\{\triv,\epsilon,\rho\}$.
By Theorem \ref{imrho},
$$
  m_\triv + m_\epsilon + m_\rho \equiv
  \ord_p \frac{C(A/F)C(A/K)^2}{C(A/M)C(A/L)^2}
  \equiv
  \ord_p \frac{C(A/F)}{C(A/M)} \mod 2.
$$
\end{example}

\begin{example}
\label{ex2}
Suppose $p$ is an odd prime and $G=\smallmatrix1*0* \subset \GL_2(\F_p)$.
Let $M, L$ be intermediate fields of degree $p-1$ and $p$ over $K$,
respectively. The group $G$ has
$p-1$ one-dimensional representations, all of which are realisable
over $\Q_p$ and factor through $\Gal(M/K)$; write $\epsilon$ for the one
of order 2. The only other irreducible representation $\rho$ of $G$
has dimension $p-1$ and can be realised over $\Q$.
There is a relation between permutation representations,
$$
  \Theta=\{1\} - \>(p\!-\!1)\Gal(F/L) - \Gal(F/M) + (p\!-\!1)\>G.
$$
Here $\cS_\Theta=\{\triv,\epsilon,\rho\}$, and Theorem \ref{imrho} implies
$$
  m_\triv + m_\epsilon + m_\rho \equiv
  \ord_p \frac{C(A/F)C(A/K)^{p-1}}{C(A/M)C(A/L)^{p-1}}
  \equiv
  \ord_p \frac{C(A/F)}{C(A/M)} \mod 2.
$$
\end{example}

\endcomment

\begin{example}
\label{ex3}
Suppose $p$ is an odd prime, and consider
$$
  B   = \smallmatrix**0*,
  U_1 = \smallmatrix\square*0*,
  U_2 = \smallmatrix**0\square
     \quad<\quad
  G=\GL_2(\F_p),
$$
where $\square$ stands for non-zero squares.
An elementary computation with
%\newpage\noindent
double cosets shows that
\beq
  B\backslash G/ U_i = B \backslash G / B &=& \{B,G \setminus B\} \cr
  U_i\backslash G/ U_i &=& \{U_i, B\setminus U_i, G \setminus B\}\cr
  U_1\backslash G/ U_2 &=& \{B, \Sigma,\Sigma' \}
\eeq
with $\Sigma$ and $\Sigma'$ the sets of matrices $\smallmatrix abcd$ with
$c(bc-ad)$ non-zero square and non-square, respectively. Since
$\langle\Ind_H^G\triv,\Ind_{H'}^G\triv\rangle=|H\backslash G/H'|$,
$$
  \Ind_B^G\triv = \triv\oplus\sigma, \quad
  \Ind_{U_1}^G\triv = \Ind_{U_2}^G\triv = \triv\oplus\sigma\oplus\rho,
$$
where $\sigma$ (Steinberg) is irreducible of dimension $p$, and
$\rho$ irreducible of \hbox{dimension} $p+1$. In particular,
$$
  \Theta = U_1-U_2
$$
is a relation between permutation representations. We will compute the
regulator constants $\RC(\Theta,\tau)$ for this relation.

If $\tau$ is not $\triv$, $\sigma$ or $\rho$, then
by Frobenius reciprocity
$$
  \dim\tau^{U_i} =
  \langle \triv_{U_i}, \Res_{U_i} \tau \rangle =
  \langle \triv+\sigma+\rho, \tau \rangle = 0,
$$
so $\RC(\Theta,\tau)=1$. On the other hand,
$$
  \RC(\Theta,\triv) = 1, \quad
  \RC(\Theta,\sigma) = 1, \quad
  \RC(\Theta,\rho) = p.
$$
To see this, it suffices to verify that
$\RC(\Theta,\triv) = \RC(\Theta,\Ind_B^G\triv) = 1$ and
$\RC(\Theta,\Ind_{U_1}^G\triv) = p$. These are permutation representations,
so they come with the standard ``identity'' pairing $\lara$. It satisfies
$$
  \det\bigl(\tfrac{1}{|N|}\lara|(\Ind_H^G\triv)^N\bigr)
    = \prod_{x\in N\backslash G/H} \tfrac{|NxH|}{|N||H|},
$$
which is easy to compute from the explicit description of double cosets.

By Theorem \ref{imrho}, for any principally polarised abelian variety $A/K$,
$$
  \rksel A{F^{U_1}}p-\rksel A{F^{B}}p =
  m_\rho \equiv
  \ord_p \frac{C(A/F^{U_1})}{C(A/F^{U_2})} \mod 2.
$$
Note that although $F^{U_1}$ and $F^{U_2}$ are arithmetically equivalent
fields (they have the same zeta-function), the right-hand side need not
be 0. For instance, take  an elliptic curve $E/\Q$
of prime conductor $l\ne p$ and split multiplicative reduction at $l$,
$p\nmid\ord_l j(E)$,
$l$ a primitive root modulo $p$, and $\Gal(\Q(E[p])/\Q)=\GL_2(\F_p)=G$.
(E.g. take $E=X_1(11)$, $p=3$). Then
the decomposition and inertia subgroups of $l$ in $G$ are
$$
  D = \smallmatrix **01 \qquad\text{and}\qquad
  I = \smallmatrix 1*01.
$$
It is easy to see that $l$ decomposes in $F^{U_1}$ and $F^{U_2}$ as
\begingroup\def\p{{\mathfrak p}}\def\q{{\mathfrak q}}%
$$
 l = \p_1 \p_2^p \p_3^p
 \qquad\text{and}\qquad
 l = \q_1 \q_2^p \q_3
$$
with $\p_1$ and $\q_2$ of residue degree 2 and the rest of residue degree 1
over $l$. So
$$
  C(E/F^{U_1}) = p^2\cdot c_l(E/\Q)^3, \qquad C(E/F^{U_2}) = p\cdot c_l(E/\Q)^3,
$$
and $m_\rho$ is therefore odd.
\endgroup
\end{example}


\begin{thebibliography}{29}

\bibitem{BS}
B. J. Birch, N. M. Stephens,
The parity of the rank of the Mordell-Weil group,
Topology 5 (1966), 295--299.

\bibitem{BL}
S. Bosch, Q. Liu,
Rational points of the group of components of a N\'eron model,
Manuscripta Math. 98 (1999), no. 3, 275--293.

\bibitem{TV-S}
T. Dokchitser, V. Dokchitser,
On the Birch--Swinnerton-Dyer quotients modulo squares,
2006, arxiv: math.NT/0610290.

\bibitem{TV-T}
T. Dokchitser, V. Dokchitser,
Regulator constants and the parity conjecture,
2007, arxiv: 0709.2852.

\bibitem{Gre}
R. Greenberg, On the Birch and Swinnerton-Dyer conjecture,
Invent. Math. 72, \linebreak no.~2 (1983), 241--265.

\bibitem{Gre2}
R. Greenberg,
Trivial zeros of $p$-adic $L$-functions,
Contemporary Math. 165 (1994), 149--174.

\bibitem{Gro}
A. Grothendieck, Mod\`eles de N\'eron et monodromie, LNM 288,
S\'eminaire de G\'eom\'etrie 7, Expos\'e IX, Springer-Verlag, 1973.

\bibitem{Guo}
L. Guo,
General Selmer groups and critical values of Hecke L-functions,
Math. Ann. 297 no. 2 (1993), 221--233.

\bibitem{Kim}
B. D. Kim,
The Parity Theorem of Elliptic Curves at Primes with
Supersingular Reduction, Compositio Math. 143 (2007) 47--72.

\bibitem{Mac}
S. MacLane, Homology, Springer-Verlag, Berlin Heidelberg New York, 1995.
(Reprint of the 1975 ed.)

\bibitem{MR}
B. Mazur, K. Rubin, Finding large Selmer ranks via an arithmetic theory
of local constants, arxiv: math.NT/0512085, to appear in Annals of Math.

\bibitem{MilO}
J. S. Milne, On the arithmetic of abelian varieties,
Invent. Math. 17 (1972), 177--190.

\bibitem{MilA}
J. S. Milne, Arithmetic duality theorems,
Perspectives in Mathematics, No. 1, Academic Press, 1986.

\bibitem{Mon}
P. Monsky, Generalizing the Birch--Stephens theorem. I: Modular curves,
Math. Z., 221 (1996), 415--420.

\bibitem{Nek}
J. Nekov\'a\v r, Selmer complexes, Ast\'erisque 310 (2006).

\bibitem{Sil2}
J. H. Silverman, Advanced Topics in the Arithmetic of Elliptic Curves,
GTM 151, Springer-Verlag 1994.




\bibitem{TatN}
J. Tate, Number theoretic background, in: Automorphic forms,
representations and L-functions, Part 2
(ed. A. Borel and W. Casselman), Proc. Symp. in Pure Math.
33 (AMS, Providence, RI, 1979) 3-26.


\end{thebibliography}
\end{document}